\documentclass[12pt]{article}
\usepackage{a4,amsmath,amsthm,amssymb,amsfonts,enumerate,hyperref}

\theoremstyle{plain}

\newtheorem{thm}{Theorem}      
\newtheorem{cor}{Corollary}
\newtheorem{lem}{Lemma}

\theoremstyle{definition}

\newcommand{\abs}[1]{{\left| {#1} \right|}}
\newcommand{\p}[1]{{\left( {#1} \right)}}

\author{Johan Andersson\thanks{Department of Mathematics, Stockholm University, SE-10691, Sweden. {\it johana@math.su.se}}}

\begin{document}

\title{Lower bounds in some power sum problems}

%\subjclass{Primary 11N30}

%\keywords{power sum, Turan method, explicit solutions}
\maketitle

\begin{abstract}
We prove that for $j \geq 0$ one has that
 \begin{gather} \label{star} \tag{*}
 \inf_{\abs{z_k}=1} \max_{\nu=1,\ldots,n^2+j} \abs{\sum_{k=1}^n z_k^\nu} \geq \sqrt{n + \frac {(1+j)(n-1)}{2(j+n^2)}}.
\end{gather}
 This improves upon former estimates. Our proof will use Fej\'er kernels. We also prove corresponding results for non pure power sums, and the pure power sum estimate
 \begin{gather*}
  \p{\sqrt{2-\frac 2 \alpha}-o(1)}   \sqrt n  \leq \inf_{\abs{z_k}=1} \max_{\nu=1,\ldots,\lfloor \alpha n^2 \rfloor} \abs{\sum_{k=1}^n z_k^\nu},
\end{gather*}
for constants $\alpha>1$. This improves further on \eqref{star} when $j \geq 2 n^2$.
\end{abstract}

\section{Introduction}

The power sum method of Tur\'an (see Tur\'an \cite{Turan} or Montgomery \cite{Montgomery} Chapter 5) allows us to obtain lower bounds for power sums
\begin{gather}
  \max_{\nu=N(n),\ldots,M(n)} \abs{g(\nu)}, \\ \intertext{where} 
  g(\nu)=\sum_{k=1}^n b_k z_k^\nu \\ \intertext{for $z_k$ and $b_k$ complex numbers, and $M(n)-N(n) \geq n$ a function of $n$. We will henceforth assume that the $b_k \geq 0$ are positive. In particular we are interested in the case of {\em pure} power sums ($b_k=1$) }
   S(\nu)= \sum_{k=1}^n z_k^\nu, 
\end{gather} and  the minimum norm $\min_{k} |z_k| = 1$. We will also assume that $N(n)=1$. In this case  a number of results  have been proved.
\begin{align*}
 \max_{\nu=1,\ldots,n} \abs{S(\nu)}&\geq 1,  \qquad &\text{(Tur\'an \cite{Turan2})} \\
    \max_{\nu=1,\ldots,2nm-m(m+1)+1} \abs{S(\nu)}&\geq \sqrt{m}.  \qquad (1 \leq m \leq n)  \qquad &\text{(Andersson \cite{Andersson})}  
\\ \intertext{Under the min norm it seems reasonable that the minimal systems 
$(z_1,\ldots,z_n)$ which minimize the expressions actually lies on, or very close to the unit circle. This has been difficult to prove and in fact in that case, when the $z_k$ are unimodular, Newman, Cassels and Szalay have independently  proved the stronger result}  
\max_{1 \leq \nu \leq c n} \left| \sum_{k=1}^n z_k^\nu \right|&\geq \sqrt{\frac {cn-n+ 1} {c}}. \qquad &\text{(\cite{Turan}, Theorem 7.3)}
\end{align*}

\section{One sided bounds}

We will denote
\begin{gather} \label{star3}
  g(\nu)=\sum_{k=1}^n b_k e(\theta_k \nu),
\end{gather}
where $\theta_k$ are real numbers and  $b_k >0$. We will let
\begin{gather*}
 A= g(0)=\sum_{k=1}^n b_k, \qquad \text{and} \qquad B=\sum_{k=1}^n b_k^2.
\end{gather*}
In this section we will furthermore assume that $g$ is real valued. In particular this implies that
\begin{gather} \label{hhh}
  g(\nu)=g(-\nu).
\end{gather}
We will also let
\begin{gather*}
 g^+(\nu)=\begin{cases} g(\nu), & g(\nu)>0, \\ 0, & \text{otherwise,} \end{cases}  \hskip -13pt \qquad \text{and} \hskip -13pt \qquad  g^-(\nu)=\begin{cases} g(\nu), & g(\nu)<0, \\ 0, & \text{otherwise.} 
\end{cases}
\end{gather*}
It is clear that
\begin{gather} \notag
 g(\nu)= g^+(\nu)+g^-(\nu), \\ \intertext{and}
 \abs{g(\nu)}= g^+(\nu)-g^-(\nu). \label{oj23}
\end{gather}

Our method of proof will use the Fej\'er kernel
 \begin{gather} \label{star1}
   F_{m+1}(x)=\sum_{\nu=-m}^{m}\p{1-\frac{\abs{\nu}}{m+1}} e(\nu x).
\\ \intertext{The Fej\'er kernel can be written as}
 \label{nonneg}
   F_{m+1}(x) =\frac 1 {m+1} \left( \frac{\sin \pi (m+1)x}{\sin \pi x} \right)^2,
\end{gather} 
and  is thus non negative. We will let 
\begin{gather} \label{starr2}
  \alpha= \frac 2 {m} \sum_{\substack{g^+(\nu)>0 \\ 1 \leq \nu \leq m}} \p{1-\frac{\nu}{m+1}}, \qquad \text{and}  \qquad \beta= \frac 2 {m} \sum_{\substack{g^-(\nu)<0 \\ 1 \leq \nu \leq m}} \p{1-\frac{\nu}{m+1}}.
\end{gather}
From the representation \eqref{nonneg} it follows that  $F_{m+1}(0)=m+1$. From this it is clear that
 \begin{gather} \label{abbb}
   \sum_{\nu=1}^m \p{1-\frac {\nu}{m+1}} = \frac m 2,
  \end{gather}
 and thus also
\begin{gather} \label{star2}
 \alpha+\beta \leq 1.
\end{gather}

\subsection{The first method}

\begin{lem}
  One has that
  \begin{gather*}
    \sum_{\nu=1}^{m} \p{1-\frac {\nu}{m+1}} \abs{g(\nu)}^2 \geq \frac{(m+1) B-A^2} 2.
  \end{gather*}
\end{lem}

\begin{proof} We have that
\begin{gather*}  
  \sum_{\nu=-m}^{m}  \p{1-\frac{\abs{\nu}}{m+1}}  \abs{g(\nu)}^2 = \sum_{k,l=1}^n b_k b_l F_{m+1}(\theta_{k}-\theta_{l}),  \\ 
\intertext{which by the contribution of the diagonal $k=l$, and the non negativity of  the Fej\'er kernel, eq. \eqref{nonneg} implies that}
 \sum_{\nu=-m}^{m}  \p{1-\frac{\abs{\nu}}{m+1}} \abs{g(\nu)}^2 \geq \sum_{k=1}^n b_k^2  F_{m+1}(0).
\end{gather*}
The result follows by subtracting the term $\nu=0$ and using equation \eqref{hhh}.
\end{proof}

\begin{lem} \label{lem2} Suppose that  $\abs{g(v)} \leq M$ for $\nu=1,\ldots,m$. Then
\begin{gather*}
     \sum_{\nu=1}^{m} \p{1-\frac {\nu}{m+1}} g^+(\nu) \geq \frac {B(m+1) - AM-A^2} {4M}.
\end{gather*}
\end{lem}

\begin{proof}
Since $g^+(\nu)=(\abs{g(\nu)}+g(\nu))/2$ and
 obviously $\abs{g(\nu)} \geq \abs{g(\nu)}^2/M$ when $\nu \neq 0$ we have that
 \begin{gather*}
       \sum_{\nu=1}^{m} \p{1-\frac {\nu}{m+1}} g^+(\nu) \geq   \frac 1 {2} \sum_{\nu=1}^{m}  \left( 1-\frac{\nu}{m+1} \right) \left(\frac{\abs{g(\nu)}^2} {M}  + g(\nu) \right)  .
\end{gather*}
The  first term can be estimated by Lemma 1 and gives the contribution
\begin{gather*}
 \frac {B(m+1) -A^2} {4M}.
\end{gather*}
The second term can be investigated by use of the Fej{\'e}r kernel. By the non  negativity of the Fej\'er kernel, equation \eqref{nonneg} we have that
\begin{gather*}
   \sum_{\nu=-m}^{m}  \left( 1-\frac{\abs{\nu}}{m+1} \right)  g(\nu) \geq 0
  \end{gather*}
By the fact that $g(0)=A$ and  using equation \eqref{hhh}  we find that
\begin{gather}
 \sum_{\nu=1}^{m}  \left( 1-\frac{\nu}{m+1} \right)  g(\nu) \geq - \frac A 2,
\end{gather}
which gives the remaining contribution to our Lemma.
\end{proof}
We now prove the following Theorem.

\begin{thm}
Suppose that $\abs{g(\nu)} \leq M$ for $\nu=1,\ldots,m$.  Then one has that
\begin{gather*}
  \max_{\nu=1,\ldots,m} g^+(\nu) \geq \frac {B(m+1)-AM-A^2} {2M m}.
\end{gather*}
\end{thm}
\begin{proof}
  This follows from Lemma 2 and equation \eqref{abbb}.
\end{proof}

\subsection{An improvement for large values of $m$}
As $m$ tends to infinity Theorem 1 will give us
\begin{gather*}
  \max_ {\nu=1,\ldots,m} g^+(\nu) \geq   \frac {B} {2M}-o(1),
\end{gather*}
We will prove a stronger results which allows us to obtain
\begin{gather*}
  \max_ {\nu=1,\ldots,m}  g^+(\nu) \geq  M+2-\frac{2M^2} B-o(1).
\end{gather*}
This will give sharper results when $M \asymp \sqrt B$ and for large $m$.

\begin{thm}
Suppose that $\abs{g(\nu)}\leq M$ for $\nu=1,\ldots,m$. If 
$B(m+1)-A^2-mM^2 \geq 0$ then  one  has that
\begin{gather*}
  \max_{\nu=1,\ldots,m} g^+(\nu) \geq \frac{B(m+1)-A^2}{mM}, \\ \intertext{In case $B(m+1)-A^2-m M^2 \leq 0$ one has that}
     \max_{\nu=1,\ldots,m}     g^+(\nu) \geq    M+2 \times \frac{B(m+1)-A^2-mM^2} {B(m+1)-A^2-AM}
\end{gather*}
 under the assumption that the denominator in the last fraction is positive.
\end{thm}

\begin{proof}
By equations \eqref{starr2} and \eqref{star2}
it is clear that
\begin{gather} \label{oj22}
  \frac 2 m \sum_{\nu=1}^{m}  \p{1-\frac{\nu} {m+1}} g^-( \nu) \geq    - M (1-\alpha).
\end{gather}
By equation \eqref{oj23} and the fact that $\abs{g(\nu)} \leq M$ we get the inequality 
\begin{gather*}
  g^+(\nu) \geq \frac{\abs{g(\nu)}^2} M+g^-(\nu).
\end{gather*} By combining this with equation \eqref{oj22} we see that
\begin{align*}
  \frac 2 m \sum_{\nu=1}^{m}  \p{1-\frac{\nu} {m+1}} g^+( \nu) &\geq \frac 2 m  \sum_{\nu=1}^{m}   \p{1-\frac{\nu} {m+1}} \frac{\abs{g( \nu)}^2} M -   M (1-\alpha).
\\ \intertext{which by Lemma 1 can be estimated by}
&\geq  \frac{B (m+1)- A^2} {Mm} -  M (1-\alpha).  \end{align*}
This together with the definition of $\alpha$, equation \eqref{starr2}  implies  that
\begin{gather} \label{ttt}
 g^+(\nu) \geq  
 \frac{1} {\alpha} \cdot \p{ \frac{B (m+1)- A^2} {Mm} -  M (1-\alpha)} =
  \frac {B(m+1)-A^2-mM^2} {\alpha M m}+M
\end{gather}
for some $\nu=1,\ldots,m$.
We see that if $B(m+1)-A^2-mM^2 \geq 0$ then the function is decreasing in $\alpha$ and the  minimum over $0 < \alpha \leq  1$ is attained for $\alpha=1$. This gives us case 1. In the case when $B(m+1)-A^2-mM^2 \leq 0$ the function is increasing in $\alpha$ and we use the following estimate
\begin{gather*}
 \alpha  \geq \frac {B(m+1)-AM-A^2}{2 m M}.
\end{gather*}
which follows from Lemma 2 to obtain a lower bound. By putting this value in the right hand side of \eqref{ttt} we obtain a lower bound and we obtain the second part of our theorem. We remark that we also need that $B(m+1)-AM-A^2$ is positive since otherwise we get an $\alpha<0$.
\end{proof}

\section{A lower bound for power sums}
We will now use our one sided theorems to obtain improved lower bounds for the absolute values of power sums.

\begin{thm}
 Let 
 \begin{gather*}
   B_\nu=\sum_{k=1}^n b_k^\nu, \qquad A=B_1^2-B_2, \qquad \text{and} \qquad B=B_2^2-B_4.
\end{gather*}
 Then one has  that 
 \begin{gather} \label{i0} \max_{\nu=1,\ldots,m} \abs{g(\nu)}^2 \geq B_2+\frac{B (1+1/m)}{2B_2}- \frac {AB_2+A^2} {2B_2 m}. \\ \intertext{One also has that}
\label{i1}\max_{\nu=1,\ldots,m}  \abs{g(\nu)}^2  \geq  2B_2 -2 \times \frac{A^2-B+ m B_4}{B(m+1)-AB_2-A^2}
\end{gather}
 when $m \geq (B-A^2)/B_4$ and  both the numerator and the denominator in the last fraction is positive (this is true for $m$ sufficiently large).
\end{thm}

\begin{proof}
 Let $g(\nu)$ be defined by equation \eqref{star3}. Then
\begin{gather*} \begin{split}
  \abs{g(\nu)}^2 &= \sum_{k=1}^n b_k^2+ \sum_{k=1}^{n^2-n} c_k e(\lambda_k \nu), \\ &= B_2+h(\nu) \end{split}
\end{gather*}
where $c_k=b_ib_j$ and $\lambda_k=\theta_i-\theta_j$ for $i \neq j$. It is clear that $h(\nu)$ is real valued and hence we can use the methods of section 2.

Let us now assume that
  $|h(\nu)| \leq B_2$ for $\nu=1,\ldots,m$. By using Theorem 1 with $M=B_2$ we have that there exist a $\nu$ with $\nu=1,\ldots,m$ such that
\begin{gather*} 
  g^+(\nu) \geq \frac {B(m+1)-AB_2-A^2} {2B_2 m}.
\end{gather*}
 This implies \eqref{i0} in case $|h(\nu)| \leq B_2$.

We have by the definition of $B$ that
\begin{gather*}
  B(m+1)-A^2-mB_2^2  =B-A^2-m B_4
\end{gather*}
which is negative if $m \geq (B-A^2)/B_4$, and it follows from Theorem 2 with $M=B_2$ that
\begin{gather*}
 h(\nu) \geq   B_2 -\frac{2A^2-2B+2mB_4}{B(m+1)-AB_2-A^2}
\end{gather*}
 for some $\nu=1,\ldots,m$. This implies \eqref{i1} in case $|h(\nu)| \leq B_2$.

 Let us now assume that $|h(\nu)| > B_2$ for some $\nu=1,\ldots,m$. 
Since $\abs{g(\nu)}^2=B_2+h(\nu) \geq 0$ this means that  $h(\nu) > B_2$ and $\abs{g(\nu)}^2 \geq 2 B_2$. We see that this implies \eqref{i0} since $2 B_2 \geq B_2+(1+1/m)B/(2B_2)$ and the third term on the right hand side in \eqref{i0} is negative. Likewise it implies \eqref{i1} since the last  term on the right hand side in \eqref{i1} is negative.
\end{proof}

\section{The pure power sum case}

In the pure power sum case we have that $B_k=n$, $A=B=n^2-n$ in Theorem 3 and it follows that

\begin{cor} One has that
\begin{enumerate}[(i)]
  \item 
  $ \displaystyle \max_{\nu=1,\ldots,m} \abs{S(\nu)}^2 \geq n+ \frac {(-1+n)(1+m-n^2)} {2 m},$ 
\item $\displaystyle
\max_{\nu=1,\ldots,m} \abs{S(\nu)}^2  \geq 2n -\frac{2(1+m-2n^2+n^3)}{(-1+n)(1+m-n^2)}. \qquad (m>n^2)$
\end{enumerate}
\end{cor}
Corollary 1 $(i)$ improves upon known results for $m$  bigger than $n^2.$ In fact it is convenient to write $m=n^2+j$ and we obtain

\begin{thm}
  One has for $j \geq 0$ that 
 \begin{gather*}
 \max_{\nu=1,\ldots,n^2+j} \abs{S(\nu)} \geq \sqrt{n + \frac {(1+j)(n-1)}{2(j+n^2)}}.
\end{gather*}
\end{thm}
For $j=0$ and in the case of unimodular numbers $z_k$ it improves slightly on the general lower bound 
\begin{gather*}
 \max_{\nu=1,\ldots,n^2} \abs{S(\nu)} \geq \sqrt{n}
\end{gather*}
in Tur\'an's problem 10 from Andersson \cite{Andersson}. We obtain

\begin{cor} One has that
\begin{gather*}
 \inf_{\abs{z_k}=1} \max_{\nu=1,\ldots,n^2} \abs{S(\nu)} \geq \sqrt{n+\frac 1 {2n}-\frac 1 {2 n^2}}.
\end{gather*}
\end{cor}
Another result where the lower bound follows from Theorem 4 and the  upper bound follows from Montgomery's construction (see Montgomery \cite{Montgomery} page 101. Example 6.) is the following
\begin{cor} Suppose that $n+1$ is a prime number. One then has that
\begin{gather*}
   \sqrt{n  + \frac 1 2 -\frac {2n-1} {2(n^2+n-1)}}\leq  \inf_{\abs{z_k}=1} \max_{\nu=1,\ldots,n^2+n-1} \abs{S(\nu)} \leq \sqrt{n+1}.
\end{gather*}
\end{cor}
We remark that the lower bound holds in general. We see that this approximately  halves the previous gap between the upper and lower bound. For further discussions of explicit constructions that yields similar upper bounds in power sum problems, see our paper Andersson \cite{Andersson3}.

In our paper Andersson \cite{Andersson2} page 17, we considered functions $\Lambda$ that fulfills 
\begin{gather*}
  \sqrt n \p{\Lambda(\alpha)-o(1)} \leq \inf_{\abs{z_k}=1}  \max_{\nu=1,\ldots,\lfloor \alpha n^2 \rfloor} \abs{S(\nu)}.
\end{gather*}
We proved that we can choose $\Lambda(\alpha)=1$ for $\alpha>0$ and furthermore that for $0<\alpha \leq 1$ we proved that it is the best possible. We asked whether the function must be identically $1$ or must be  bounded. While we can not answer if there exist such an  unbounded function it follows from Corollary 1 that we can choose a $\Lambda$ such that $\lim_{\alpha \to \infty} \Lambda(\alpha)=\sqrt 2$. More specifically we obtain the following Theorem.

\begin{thm} Let $\alpha \geq  1$ be a constant. One then has that
\begin{gather*}
  \p{\sqrt{\Phi(\alpha)}-o(1)}   \sqrt n  \leq \inf_{\abs{z_k}=1} \max_{\nu=1,\ldots,\lfloor \alpha n^2 \rfloor} \abs{S(\nu)} \leq  \p{\sqrt{\lceil \alpha \rceil}+o(1)} \sqrt n,
\\ \intertext{where}
        \Phi(\alpha)= \begin{cases} \frac 3 2-\frac 1 {2\alpha}, & 1 \leq \alpha \leq 3, \\
                                    2-\frac 2 \alpha, & \alpha \geq 3. \end{cases}
\end{gather*}

\end{thm}

\begin{proof}
 The lower bound for $1 \leq \alpha\leq 3$   follows from Corollary 1 $(i)$ with $m = \lfloor \alpha n^2 \rfloor$.
 The lower bound for $3 \leq \alpha$   follows from Corollary 1 $(ii)$ with $m = \lfloor \alpha n^2 \rfloor$.
  The upper bound follows from Theorem 6 in Andersson \cite{Andersson2}.
\end{proof}
In particular this will give us
\begin{gather*}
  \p{\sqrt{\frac 5 4}-o(1)}   \sqrt n  \leq \inf_{\abs{z_k}=1} \max_{\nu=1,\ldots,2n^2} \abs{S(\nu)} \leq  \p{\sqrt 2 +o(1)} \sqrt n.
\end{gather*}
We see that the lower and upper bounds are not the same and we do not yet have the true asymptotic. This contrasts to the case when we take the maximum over the interval $\nu=1,\ldots,n^2$ where we proved (see Andersson \cite{Andersson2})
\begin{gather*}
  \inf_{\abs{z_k}=1} \max_{\nu=1,\ldots,n^2} \abs{S(\nu)} \sim  \sqrt n.
\end{gather*}

\bibliographystyle{alpha}

\end{document}